\documentclass[preprint,12pt,3p]{elsarticle}

\usepackage{graphicx,graphics,amssymb,epsfig,amsmath,caption,subcaption}

\begin{document}

\begin{frontmatter}



\title{A regularization approach for solving Poisson's equation with singular charge sources and diffuse interfaces}
\author[label1]{Siwen Wang}
\address[label1]{Department of Mathematics, University of Alabama, Tuscaloosa, AL 35487,USA}
\author[label1]{Arum Lee}
\author[label2]{Emil Alexov}
\address[label2]{Department of Physics and Astronomy, Clemson University, Clemson, SC 29634, USA}
\cortext[cor1]{Corresponding author}
\author[label1]{Shan Zhao\corref{cor1}}
\ead{szhao@ua.edu}


\begin{abstract}
Singular charge sources in terms of Dirac delta functions present a well-known numerical challenge for solving Poisson's equation. For a sharp interface between inhomogeneous media, singular charges could be analytically treated by fundamental solutions or regularization methods. However, no analytical treatment is known in the literature in case of a diffuse interface of complex shape. This letter reports the first such regularization method that represents the Coulomb potential component analytically by Green's functions to account for singular charges. The other component, i.e., the reaction field potential, then satisfies a regularized Poisson equation with a smooth source and the original elliptic operator. The regularized equation can then be simply solved by any numerical method. For a spherical domain with diffuse interface, the proposed regularization method is numerically validated and compared with a semi-analytical quasi-harmonic method. 

\end{abstract}

\begin{keyword}
Poisson's equation; Singular source; Diffuse interface; Regularization; Green's function. 

\vspace{0.3in}
\end{keyword}
\end{frontmatter}

\section{Introduction}
The Poisson equation, as a mean field model, is widely used for the study of electrostatic interactions in biological and chemical systems at molecular level \cite{Honig95} and also for the design of semiconductor devices at the nanoscale \cite{Taur}. In typical applications, two dielectric materials are concerned in the system and one of them carries fixed point charges, which are represented as Dirac delta functions in the source term of Poisson's equation. In classical settings, a sharp  interface is assumed to separate two media, which yields a piecewise constant for the dielectric coefficient of Poisson's equation. 

Recently, the use of diffuse interface Poisson models becomes popular \cite{Li13,Xue17}. For biological and chemical systems in molecular or nano scales,  the assumption of a sharp interface as the boundary of two dielectric materials seems to be unphysical \cite{Hazra19}. The diffuse interface model \cite{Teigen09}, in which a smooth transition layer is assumed at material boundaries, provides an alternative to model the dielectric function. For example, in studying charged objects immersed in liquids, a smooth implicit solvent model has been developed by incorporating the structures of water dipoles and ions into mean field modeling, so that the effective dielectric coefficient becomes a smoothly variant function \cite{Abrashkin07}. In studying electrostatic interactions of macromolecule and solvent, various free energy variational models have been proposed, including minimal molecular surface \cite{Bates08}, level set \cite{Cheng07}, and field phase \cite{Zhao13}. These models all feature a diffuse interface type dielectric boundary. In this letter, a simple Poisson's equation involving inhomogeneous media and a diffuse interface is studied, without invoking additional features of the above mentioned physical models. In particular, we will assume constant dielectric values inside each dielectric medium, while the dielectric function varies smoothly from one medium to another, through a narrow transition band. 

The accurate treatment of singular charge sources of Poisson's equation is a well-known challenge. Mathematically, the fixed point charges are expressed in terms of the Dirac delta functions, which are unbounded at charge centers. In conventional numerical algorithms, a trilinear scheme is often used to distribute point charges to their neighboring grid points. This is known to be a very poor approximation, and motivates a recent development of a second order accurate geometric discretization of the multidimensional Dirac delta distribution \cite{Egan17}. We note that the numerical difficulty for representing singular functions via discrete finite values could be completely avoided if charge singularities are treated \emph{analytically}. 

For Poisson's equation with singular charges and diffuse interfaces, a family of semi-analytical methods have been proposed \cite{Xue17} for eleven orthogonal coordinate systems in which the three-dimensional (3D) Laplace equation is separable. The dielectric function is assumed to be variant only in one orthogonal direction, and the underlying diffuse interface can then be approximated via several pieces of quasi-harmonic diffuse interfaces. For each quasi-harmonic dielectric function, Green's functions for Poisson's equation can be calculated analytically. The singular charges are treated analytically in this approach with diffuse interfaces. Nevertheless, this semi-analytical method is limited to simple geometries. No analytical procedure is available in the literature for singular charges with complex domains and diffuse interfaces.

In a related field, a series of regularization methods have been developed for solving the Poisson-Boltzmann (PB) equation with singular charges and sharp interfaces over any domain, see the references in \cite{Geng17}. In regularization methods, the potential function is decomposed into a singular component plus one or two other components. Satisfying a Poisson equation with the same singular sources, the singular component can be analytically solved as Coulomb potentials or Green's functions. After removing the singular part, the other potential components are bounded, and thus can be accurately solved by finite difference or finite element methods. However, all existing regularization methods are designed for piecewise constant dielectric functions with sharp interfaces. It is unclear if regularization formulation could be established for diffuse interfaces. 

This letter presents the first regularization method in the literature that is able to handle diffuse interfaces. Besides a decomposition of potential function, the success of the new method lies in a decomposition of the inhomogeneous dielectric function. The singular charge sources containing in a complex domain can then be analytically treated. The details of the proposed regularization formulation will be discussed in Section \ref{sec.theory}. This new method can be combined with any numerical discretization, and is expected to find extensive applications for various real world problems. Numerical validation for a simple example will be considered in Section \ref{sec.numerical}. Finally, this letter ends with a conclusion.

\section{Regularization formulation}\label{sec.theory}
\begin{figure*}[!tb]
	\begin{center}
		\begin{tabular}{cc}
			\includegraphics[width=0.4\textwidth]{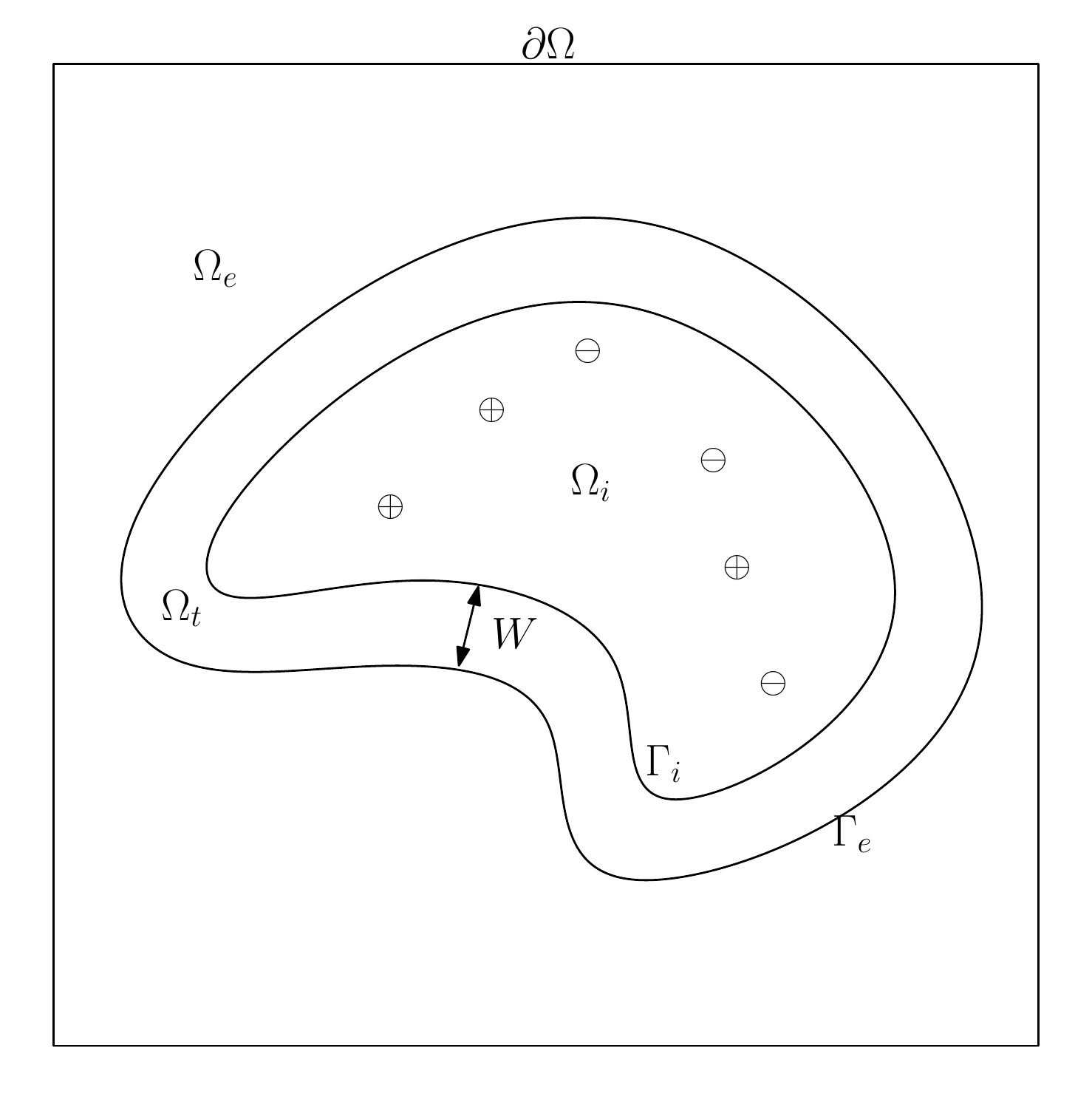} & \quad 
			\includegraphics[width=0.4\textwidth]{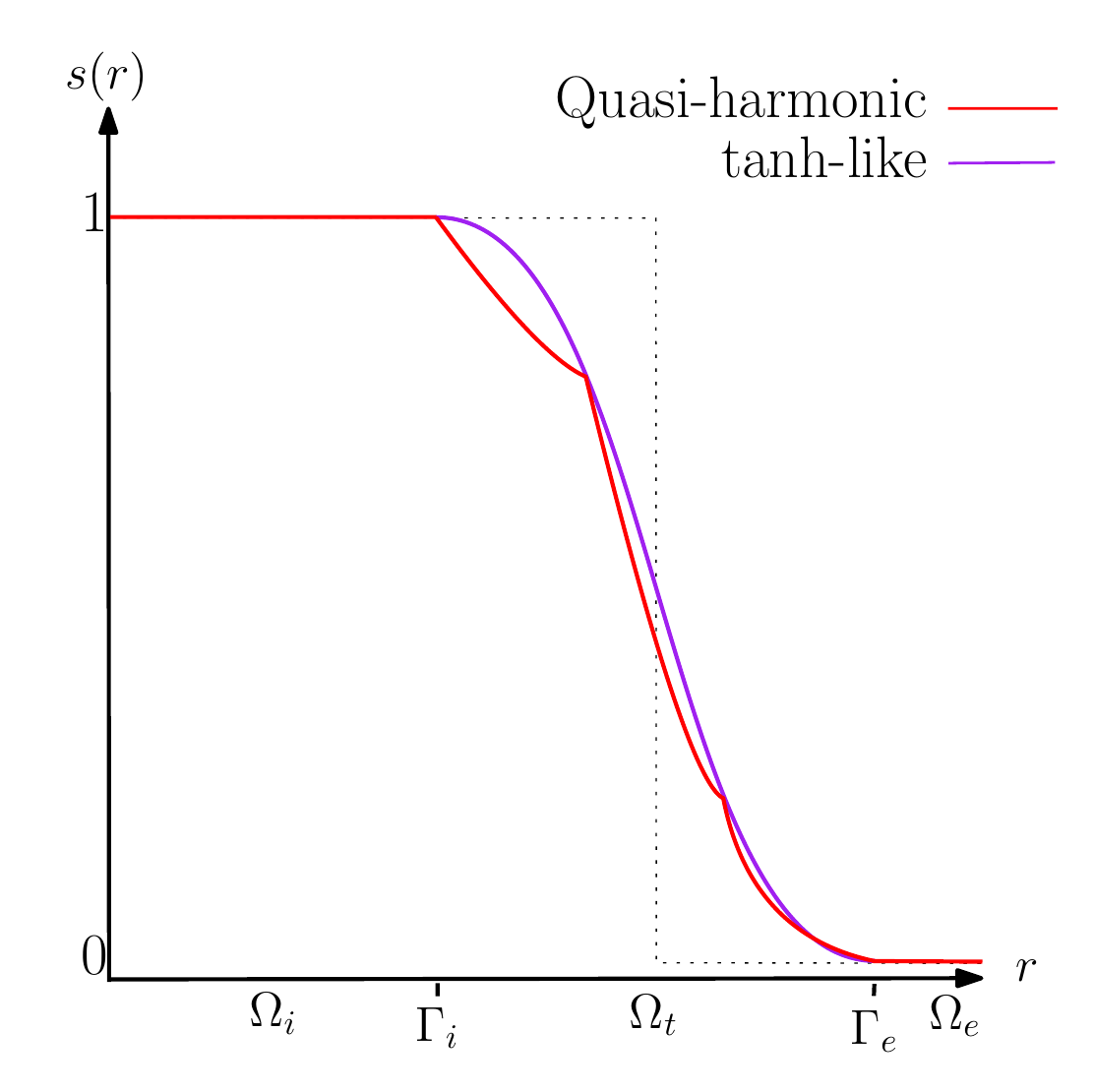} \\
			(a) & \quad  (b) 
		\end{tabular}
	\end{center}
	\caption{(a) Domain setting of the problem; (b) Diffuse interfaces used in the spherical domain example. In the regularization method, a level set function is analytically constructed via a $\tanh(\cdot)$ function, and is illustrated along the radial direction $r$. Obviously, such tanh-like diffuse interface will produce a smooth dielectric function. In the semi-analytical approach, the tanh-like diffuse interface is approximated by three pieces of quasi-harmonic diffuse interfaces.}
	\label{fig_model}
\end{figure*}

Consider a three-dimensional (3D) Poisson's equation with a Dirichlet boundary condition \cite{Xue17}
\begin{equation}\label{PE}
\begin{cases}\displaystyle
-\nabla\cdot(\epsilon(\vec{r})\nabla u(\vec{r}))= \rho :=
4\pi \sum_{j=1}^{N_s} q_j \delta (\vec{r} - \vec{r}_j), & \text{in $\Omega$},\\
u (\vec{r})= g(\vec{r}) & \text{on }\partial\Omega,
\end{cases}
\end{equation}
where $u$ is the potential function and $g$ is a boundary function. 
The domain $\Omega$ consists of three regions, an interior domain $\Omega_i$, an exterior domain $\Omega_e$, and a transition layer $\Omega_t$ in between $\Omega_i$ and $\Omega_e$. See Fig. \ref{fig_model} (a).  The interface between $\Omega_i$ and $\Omega_t$ is denoted by $\Gamma_i$, while the one between $\Omega_t$ and $\Omega_e$ is $\Gamma_e$. There exist $N_s$ point charges inside $\Omega_i$ with charge numbers being $q_j$, for $j=1,2, \ldots,N_s$.  The dielectric function $\epsilon(\vec{r})$ takes constant values $\epsilon  = \epsilon_i$ in $\Omega_i$ and $\epsilon  = \epsilon_e$ in $\Omega_e$. Here we assume $\epsilon_i < \epsilon_e$. In $\Omega_t$, $\epsilon(\vec{r})$ varies smoothly from $\epsilon_i$ to $\epsilon_e$, so that it is a $C^2$ continuous function over the entire domain $\Omega$. Consequently, function $u$ and its gradient $\nabla u$ are continuous everywhere in $\Omega$, except at charge centers. 

In the proposed two-component regularization, the potential $u$ is decomposed into a Coulomb component $u_C$ and a reaction field component $u_{RF}$ with $u=u_C + u_{RF}$. As in the sharp interface case \cite{Geng17}, the Coulomb potential is assumed to satisfy a homogeneous Poisson's equation with the same singular charges $\rho$
\begin{equation} \label{uC_eq}
\left\{
\begin{array}{ lll }
-\epsilon_i \Delta u_C({\vec r})&=\rho({\vec r}) & \hbox{in $\mathbb{R}^3$;} \\
u_C({\vec r})
&=0. & \hbox{as~}|{\vec r}| \to \infty. \\
\end{array} \right.
\end{equation}
Thus, the singular component $u_C$ is analytically given as the Green's function $G(\vec{r})$
\begin{equation}\label{uC}
u_C({\vec r})= G({\vec r}) := \sum^{N_s}_{j=1}\frac{q_j}{\epsilon_i |\vec{r} - \vec{r}_j |}.
\end{equation}

To deal with the diffuse interface, we propose to decompose the dielectric function into a constant base value plus a variant part, i.e., $\epsilon = \epsilon_i + \hat{\epsilon}$. Consequently, $\hat{\epsilon} = 0$ in $\Omega_i$ and $\hat{\epsilon} = \epsilon_e - \epsilon_i$ in $\Omega_e$, with $\hat{\epsilon} \ge 0$ throughout the domain $\Omega$. By introducing the dual decomposition into Poisson's equation (\ref{PE}), we have
\begin{equation}\label{4parts}
 -\nabla\cdot(\hat{\epsilon}\nabla u_C)-\nabla\cdot(\hat{\epsilon}\nabla u_{RF})
-\epsilon_i\Delta u_C-\epsilon_i\Delta u_{RF} =\rho, \quad \mbox{in}~\Omega.
\end{equation}
By subtracting (\ref{uC_eq}) from (\ref{4parts}), the Poisson equation is now free of singular sources
\begin{equation}\label{3parts}
 -\nabla\cdot(\hat{\epsilon}\nabla G)-\nabla\cdot(\hat{\epsilon}\nabla u_{RF})
-\epsilon_i\Delta u_{RF} =0, \quad \mbox{in}~\Omega,
\end{equation}
where we have substituted $u_C$ by the known Green's function $G$. Note that $G(\vec{r})$ is unbounded at charge centers inside $\Omega_i$. However, in  the proposed regularization, we have deliberately designed a nice property: $\hat{\epsilon} = 0$ in $\Omega_i$. This enables us to simplify the new source term of Eq. (\ref{3parts}) as,
\begin{equation}\label{new-source}
\nabla\cdot(\hat{\epsilon}\nabla G) = \nabla\hat{\epsilon}\cdot\nabla G+\hat{\epsilon}\Delta G = \nabla\hat{\epsilon}\cdot\nabla G = \nabla \epsilon \cdot\nabla G.
\end{equation}
In Eq. (\ref{new-source}), $\hat{\epsilon}\Delta G$ is dropped out, because  $\Delta G = 0$ everywhere except at charge centers within $\Omega_i$, while $\hat{\epsilon}=0$ in $\Omega_i$. In the last step, we have $\nabla \epsilon= \nabla\hat{\epsilon}$ because $\epsilon$ and $\hat{\epsilon}$ differs by a constant $\epsilon_i$. The gradient of Green's function is analytically given as
\begin{equation}
\nabla G(\vec{r}) = - \sum^{N_s}_{j=1}\frac{q_j (\vec{r} - \vec{r}_j)}{\epsilon_i |\vec{r} - \vec{r}_j |^3}.
\end{equation}
Moreover, by the definition of $\epsilon$, $\nabla \epsilon$ is non-vanishing only in $\Omega_t$, while $\nabla \epsilon=0$ for both $\Omega_i$ and $\Omega_e$.  Thus, $\nabla \epsilon \cdot\nabla G$ is finite in $\Omega$, and one just needs to calculate it in the transition band $\Omega_t$. 

In summary, we propose a new regularized Poisson's equation for the reaction field potential 
\begin{equation}\label{RPE}
\begin{cases}\displaystyle
-\nabla\cdot(\epsilon(\vec{r})\nabla u_{RF}(\vec{r}))= \nabla \epsilon(\vec{r}) \cdot\nabla G(\vec{r}), & \text{in $\Omega$},\\
u_{RF}(\vec{r}) = g(\vec{r})-G(\vec{r}) & \text{on }\partial\Omega,
\end{cases}
\end{equation}
in which the two $u_{RF}$ terms in Eq. (\ref{4parts}) have been combined into one. Hence, the decomposition of dielectric function $\epsilon = \epsilon_i + \hat{\epsilon}$ is used only in the derivation. All real computations can be carried out based on $\epsilon(\vec{r})$ only. Once $u_{RF}$ is computed from (\ref{RPE}), the solution of the original Poisson's equation (\ref{PE}) is recovered by $u=u_{RF}+G$. 

\section{Numerical validation}\label{sec.numerical}
In this letter, we validate the proposed regularization method by considering a simple geometry, i.e., a sphere. This enables us to benchmark the new method with the semi-analytical approach developed by Xue and Deng \cite{Xue17}. Moreover, the classical trilinear method is also deployed for a comparison, in which  a singular source is distributed to the vertexes of the cube or element containing the source point via a trilinear approximation. Like trilinear method, the proposed regularization method can easily handle complex geometries - such study is in progress and will be reported elsewhere. 

Consider a spherical domain $\Omega_i$ with a point charge $q_1$ at its center. Assume the charge point to be the origin of our coordinate, i.e., $\vec{r}_1=(0,0,0)$. In this example, both boundaries $\Gamma_i$ and $\Gamma_e$ are spheres with radii being $r_i=2$ and $r_e=5$, respectively. A cubic computational domain $\Omega=[-10,10]^3$ is employed. In the present study, the diffuse interface is  constructed through a level set function $s(\vec{r})$
\begin{equation}
s(\vec{r})=
\begin{cases}
s_i,& \text{if $|\vec{r}| \leq r_i$},\\
\displaystyle (s_e - s_i)\frac{ \tanh(k(\frac{|\vec{r}|-r_i}{r_e-r_i}-0.5)+1}{2}+s_i, 
& \text{if $ r_i < |\vec{r}| < r_e $},\\
s_e,& \text{if $|\vec{r}| \geq r_e$},
\end{cases}    
\end{equation}
where $k=6$ is large enough to ensure that $s(\vec{r})$ can be numerically assumed as a smooth function across $\Gamma_i$ and $\Gamma_e$. Here we take $s_i=1$ and $s_e=0$. An illustration of $s(r)$ for $r=|\vec{r}|$ is shown in Fig. \ref{fig_model} (b) as a tanh-like curve. The smooth dielectric function can then be calculated as $\epsilon(\vec{r}) = s(\vec{r}) \epsilon_i + (1-s(\vec{r})) \epsilon_e$, in which we take $ \epsilon_i=1$ and $\epsilon_e=80$.

Two numerical methods are considered. In the proposed method, the regularized Poisson's equation (\ref{RPE}) is discretized by using the central finite difference without worrying the source term singularity. In the trilinear method, the original Poisson's equation (\ref{PE}) is numerically solved by using the same finite difference discretization after trilinear distribution of the source term. In both methods, a uniform grid with the same mesh size in all three dimensions, i.e., $N=N_x=N_y=N_z$, is used with spacing $h=\frac{20}{N-1}$. On boundary $\partial \Omega$, the Dirichlet boundary data is given by the Coulomb potential for exterior medium, i.e., $g(\vec{r})=\frac{q_1}{\epsilon_e |\vec{r}|}$. In the regularization method, the final numerical solution consists of $u_{RF}(\vec{r})$ values over $N^3$ grid nodes. For the trilinear solution, in order to directly compare with $u_{RF}(\vec{r})$, we will subtract the potential solution by Green's function (\ref{uC}), and denote the resulted solution as $u_{TL}(\vec{r})$.

Generally speaking, the Poisson equation with a singular source and diffuse interface cannot be solved analytically. Fortunately, since the present geometry is separable, a semi-analytical method \cite{Xue17} is available to provide series solutions in the case of a quasi-harmonic diffuse interface. Following \cite{Xue17}, we first approximate the present tanh-like diffuse interface $s(r)$ by three pieces of quasi-harmonic diffuse interfaces. In particular, we will divide the transition region $\Omega_t$ into three spherical shells of the same thickness. Referring to Fig. \ref{fig_model} (b), this amounts to cut the interval $r \in [r_i,r_e]$ into three subintervals of equal length. Then in each subinterval, one approximates $s(r)$ by a quasi-harmonic function such that its endpoint values agree with $s(r)$. Note that the diffuse interface used in the quasi-harmonic method is then piecewise continuous and is an approximation of our tanh-like diffuse interface. For the present domain setting with a spherical domain $\Omega_i$, three spherical shells with quasi-harmonic diffuse interfaces, and a external domain $\Omega_e$, analytical series solution can be established for the Poisson equation \cite{Xue17}. Similar to the trilinear case, we will subtract the series solution by Green's function, and denote the resulted function as $u_{QH}(\vec{r})$. 

Before we present numerical results, it should be pointed out that the asymptotic limits of numerical solution and semi-analytical solution are different, because the quasi-harmonic diffuse interface is different from the tanh-like diffuse interface. Theoretically, two numerical solutions $u_{RF}(\vec{r})$ and $u_{TL}(\vec{r})$ should converge to the same place, as $h$ goes to zero. However, the difference between $u_{RF}(\vec{r})$ and $u_{QH}(\vec{r})$ will not become smaller for a smaller $h$, and could be reduced only if more pieces of quasi-harmonic functions are employed for diffuse interface approximation. 

\begin{figure}[t]
\centering
\begin{subfigure}{.32\textwidth}
  \centering
  \includegraphics[width=1.0\linewidth]{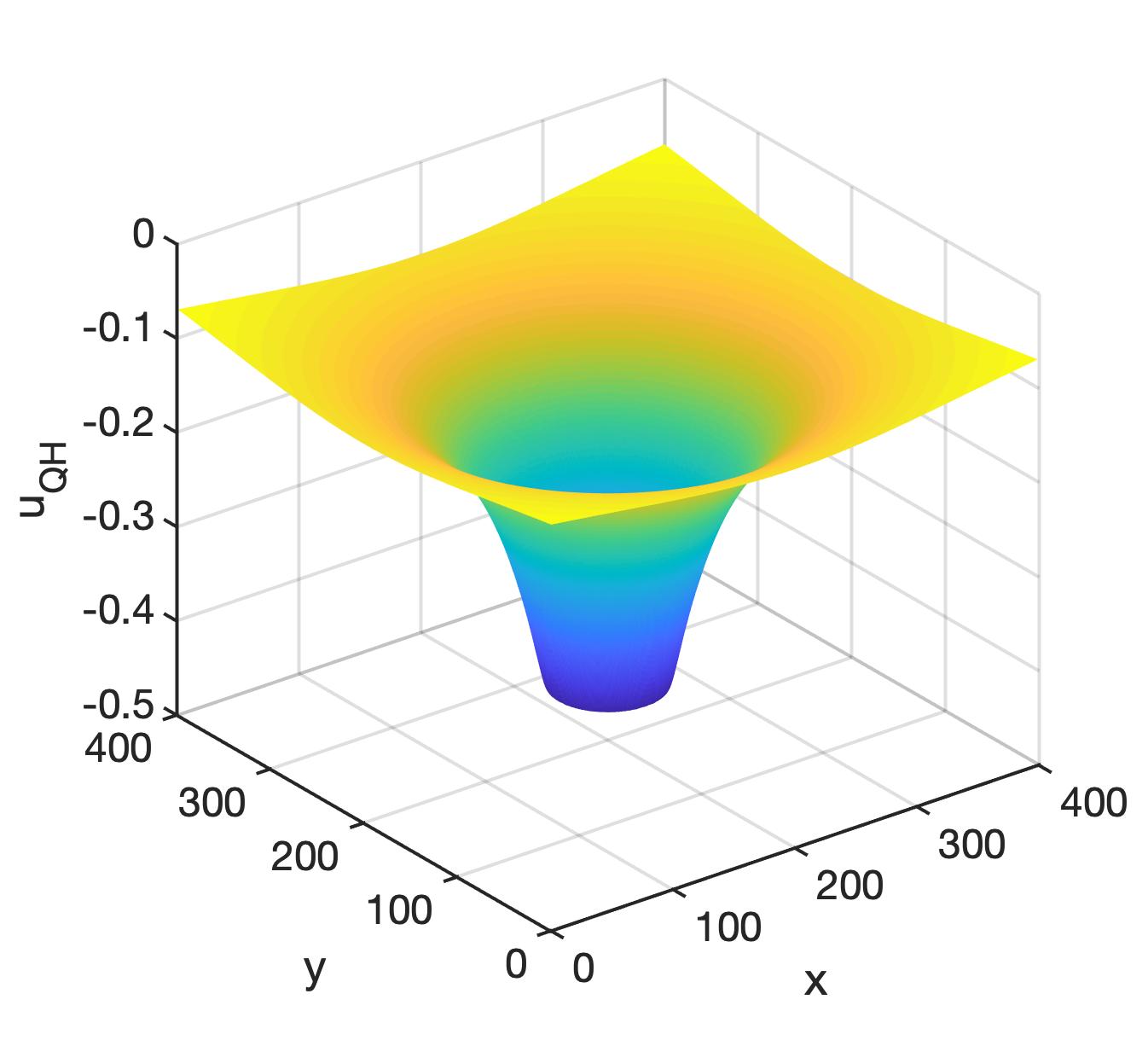}
  \caption{$u_{QH}$}
\end{subfigure}
\begin{subfigure}{.32\textwidth}
  \centering
  \includegraphics[width=1.0\linewidth]{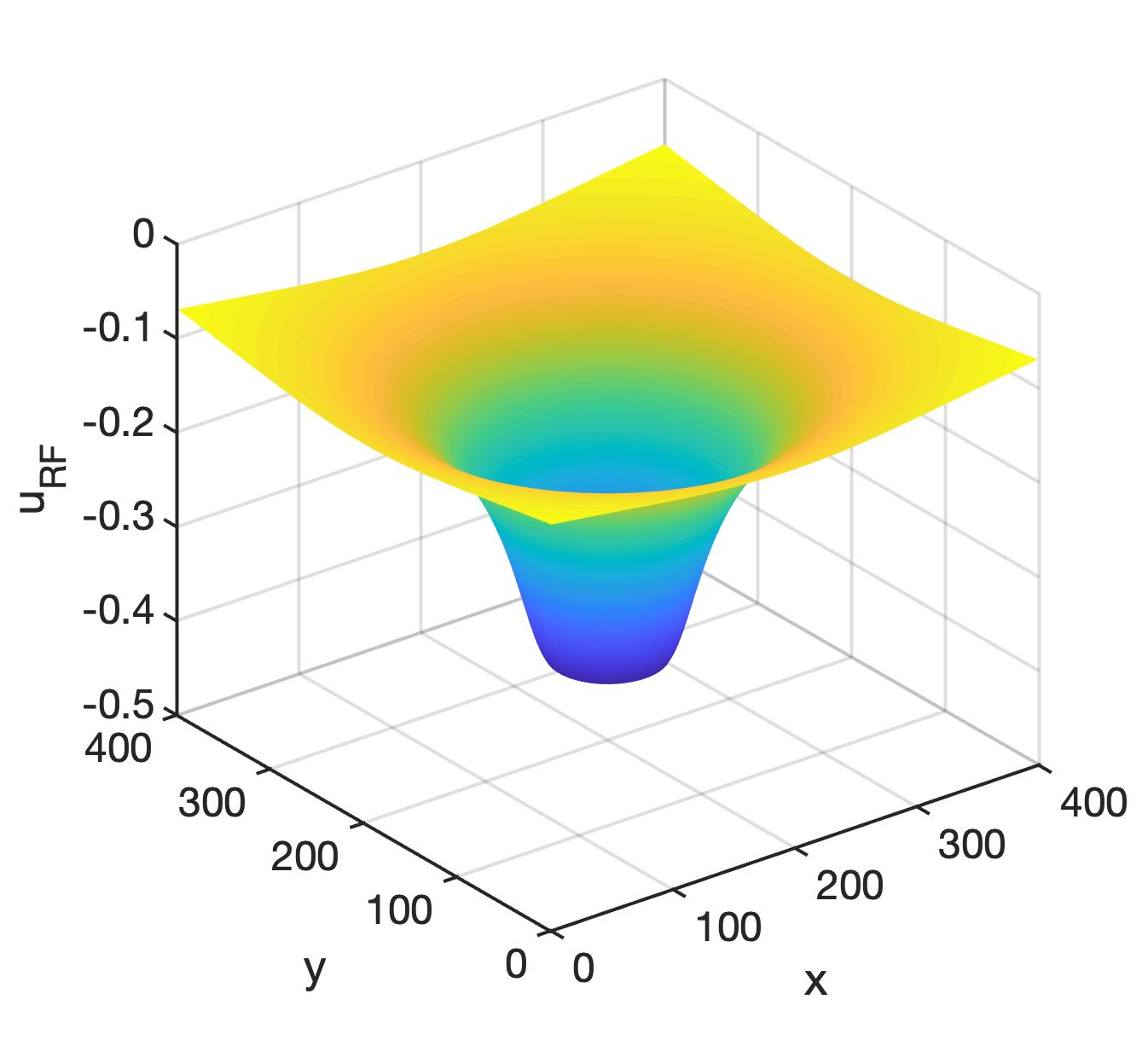}
  \caption{$u_{RF}$}
\end{subfigure}
\begin{subfigure}{.32\textwidth}
  \centering
  \includegraphics[width=1.0\linewidth]{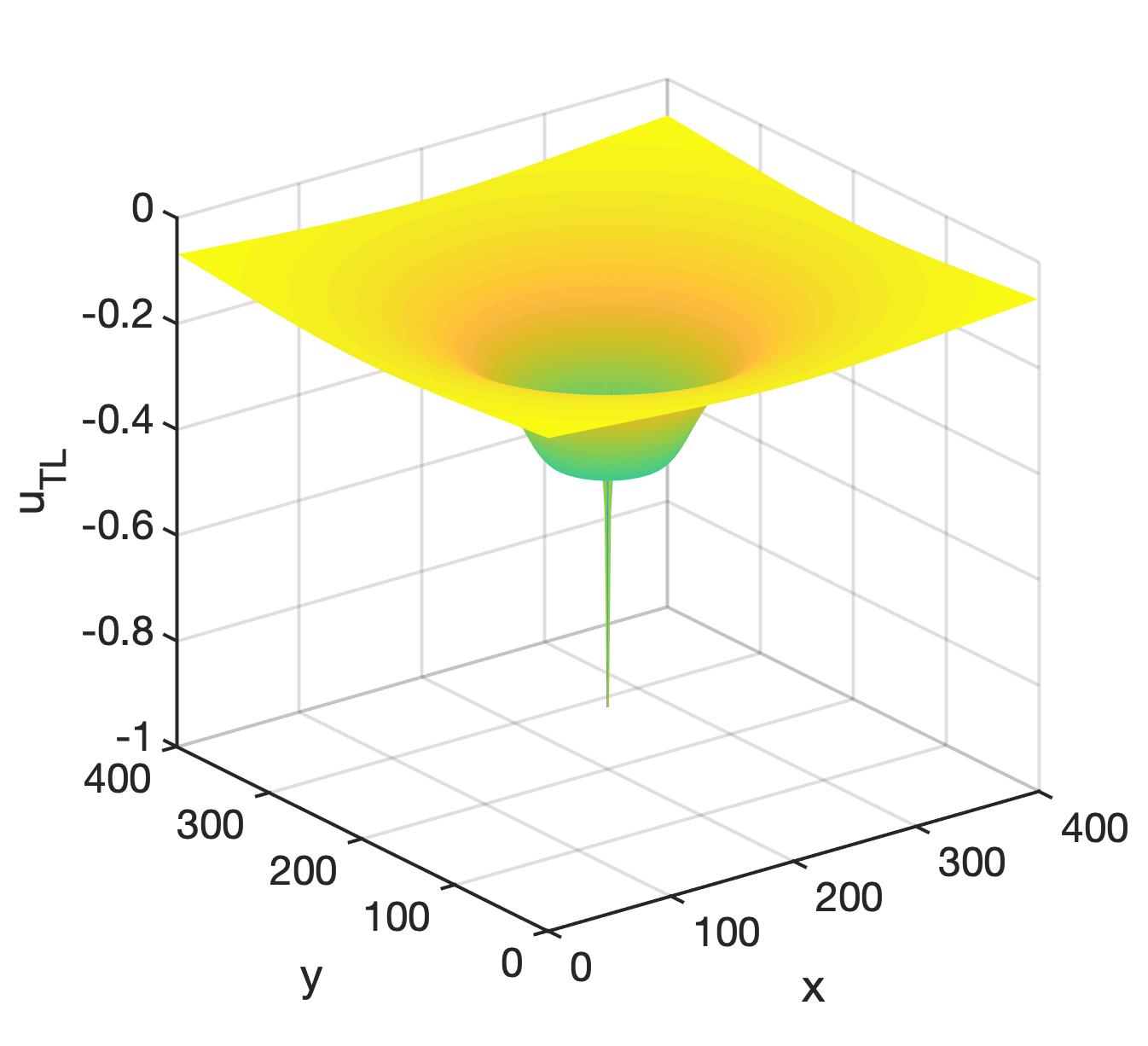}
  \caption{$u_{TL}$}
\end{subfigure}
\caption{Surface plot of potential solutions $u_{QH}$, $u_{RF}$ and $u_{TL}$ on the plane $z = 0$.}
\label{Fig.3d1}
\end{figure}

We first visually compare three solutions. By taking  $N=400$, surface plots of potential solutions $u_{QH}$, $u_{RF}$ and $u_{TL}$ on the plane $z = 0$ are shown in Fig. \ref{Fig.3d1}. It can be seen that the three solutions are almost identical for the majority part of the domain. For trilinear solution $u_{TL}$, numerical artifact is very obvious at the charge center. Excluding a small neighborhood around the origin, the difference between $u_{RF}$ and $u_{TL}$ is then very small. The regularization solution $u_{RF}$ and semi-analytical solution $u_{QH}$ have almost the same shape - a flat potential inside $\Omega_i$ with a smoothly increment outside the sphere. This indicates that the charge singularity is well taken care of in the proposed regularization method, just as in the semi-analytical quasi-harmonic method \cite{Xue17}.

\begin{figure}[!t]
\centering
\begin{subfigure}{.24\textwidth}
  \centering
  \includegraphics[width=1.0\linewidth]{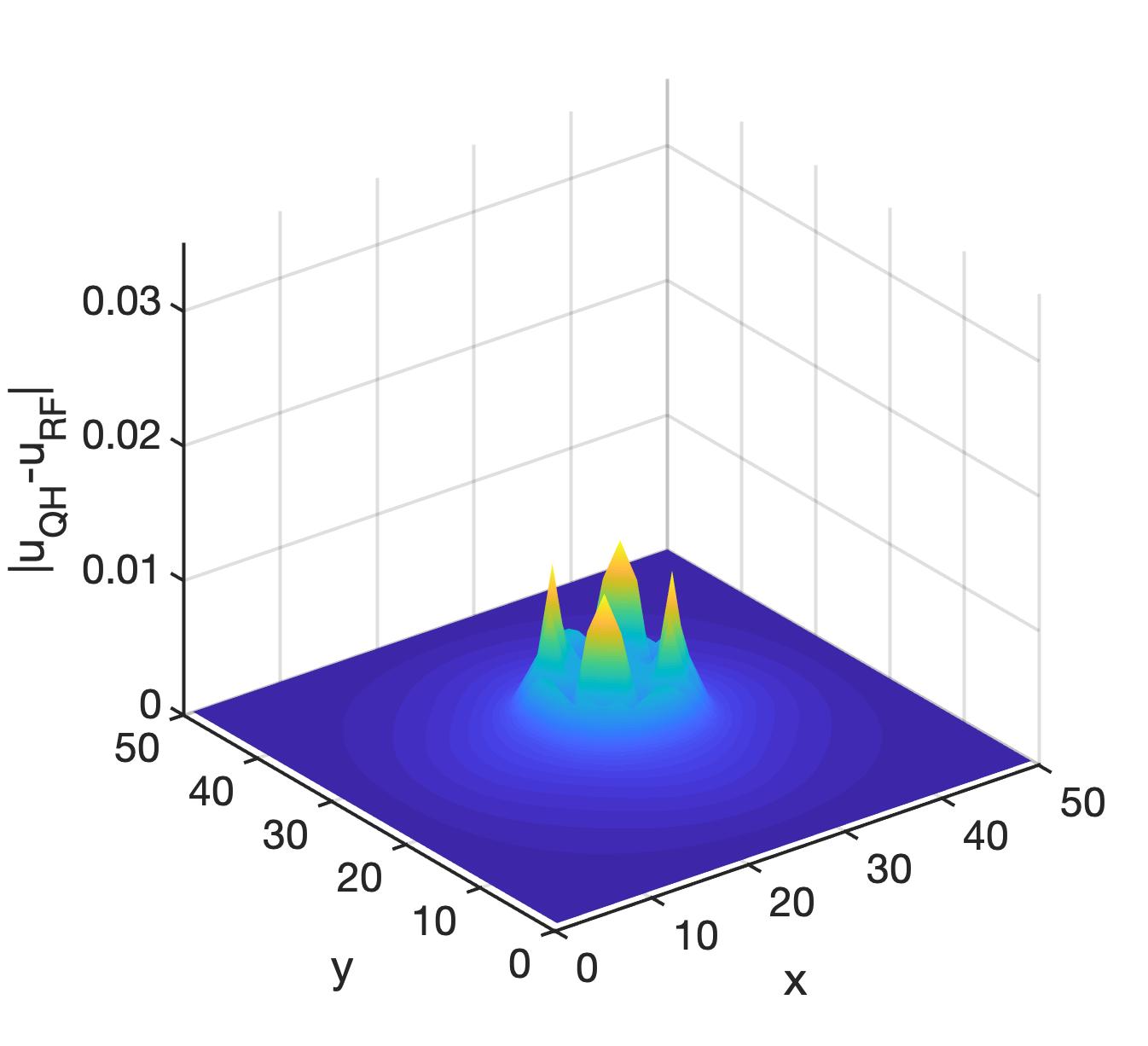}
  \caption{$N=50$}
\end{subfigure}
\begin{subfigure}{.24\textwidth}
  \centering
  \includegraphics[width=1.0\linewidth]{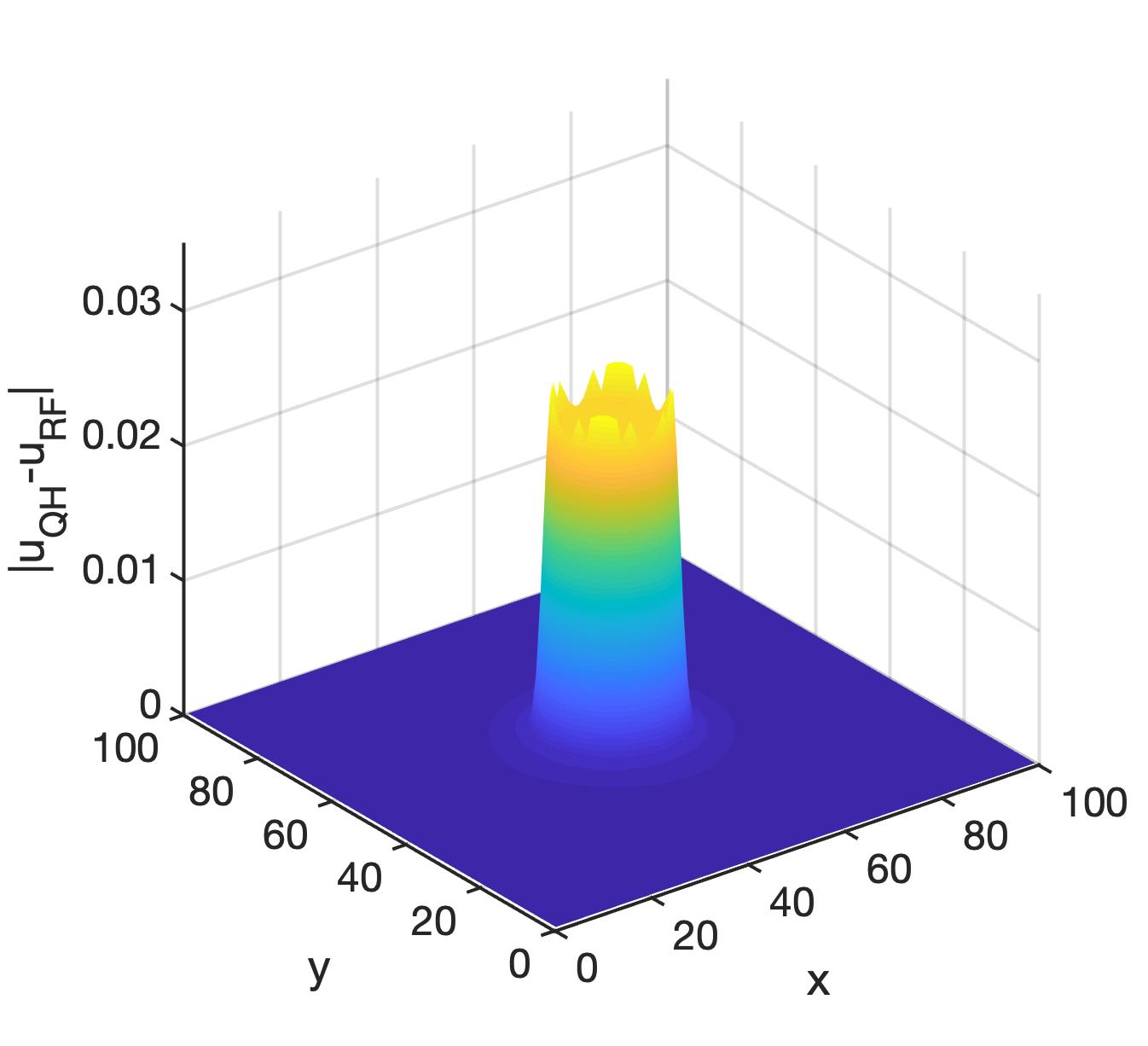}
  \caption{$N=100$}
\end{subfigure}
\begin{subfigure}{.24\textwidth}
  \centering
  \includegraphics[width=1.0\linewidth]{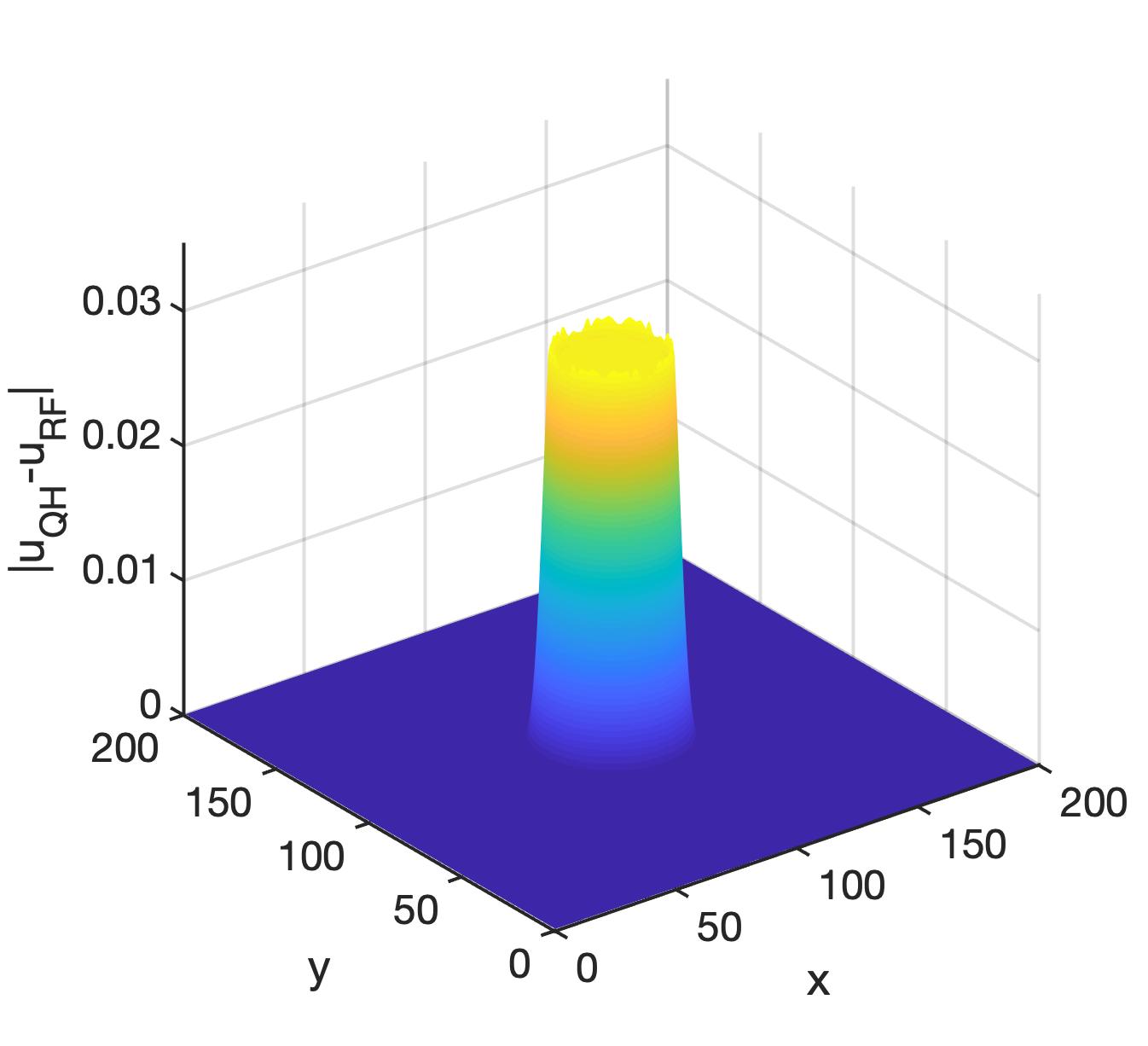}
  \caption{$N=200$}
\end{subfigure}
\begin{subfigure}{.24\textwidth}
  \centering
  \includegraphics[width=1.0\linewidth]{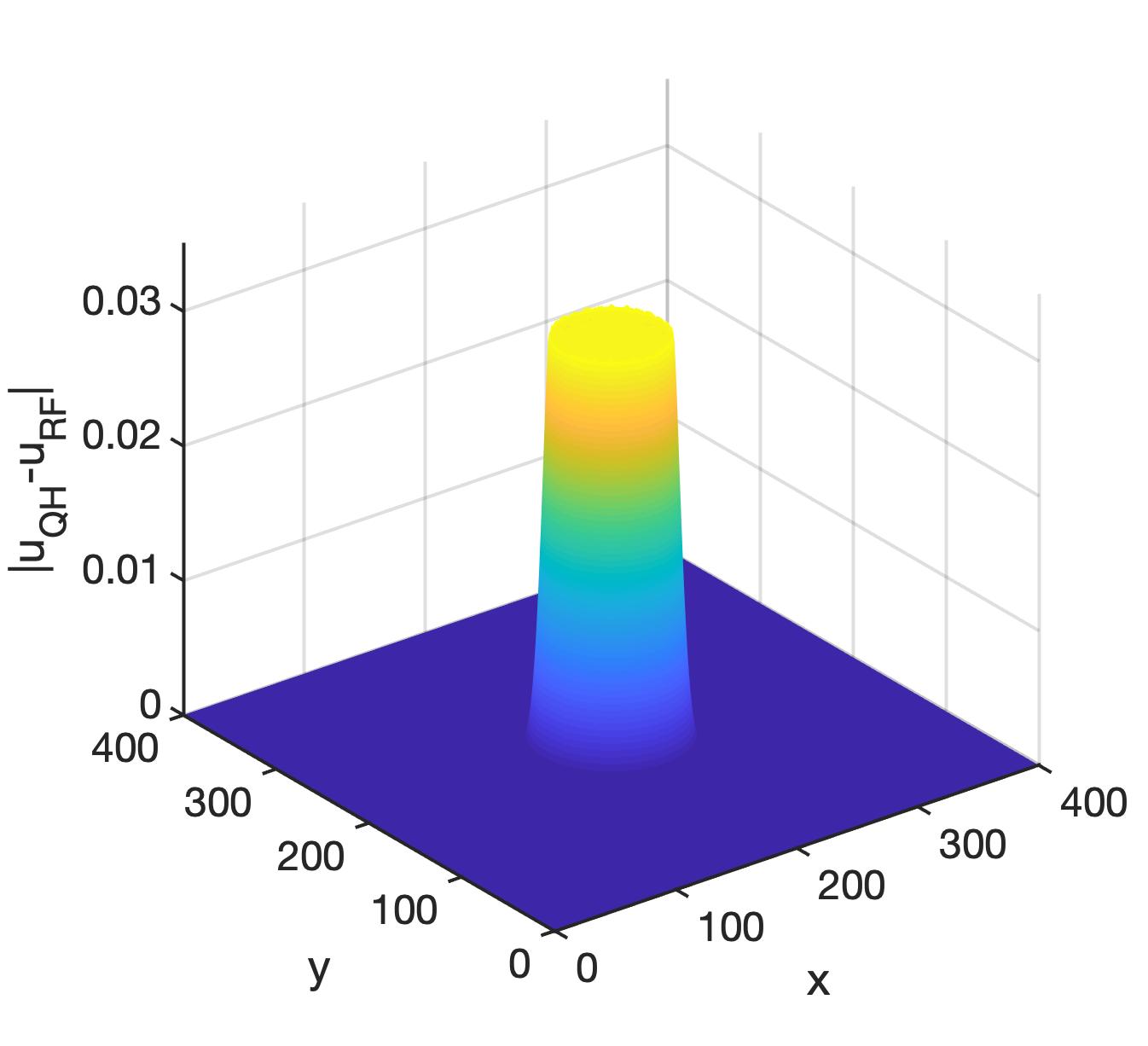}
  \caption{$N=400$}
\end{subfigure}
\caption{Surface plots of the potential difference $|u_{QH} - u_{RF}|$  on the plane z = 0 under different mesh sizes. }
\label{Fig.3d2}
\end{figure}

\begin{figure}[!t]
\centering
\begin{subfigure}{.45\textwidth}
  \centering
  \includegraphics[width=0.95\linewidth]{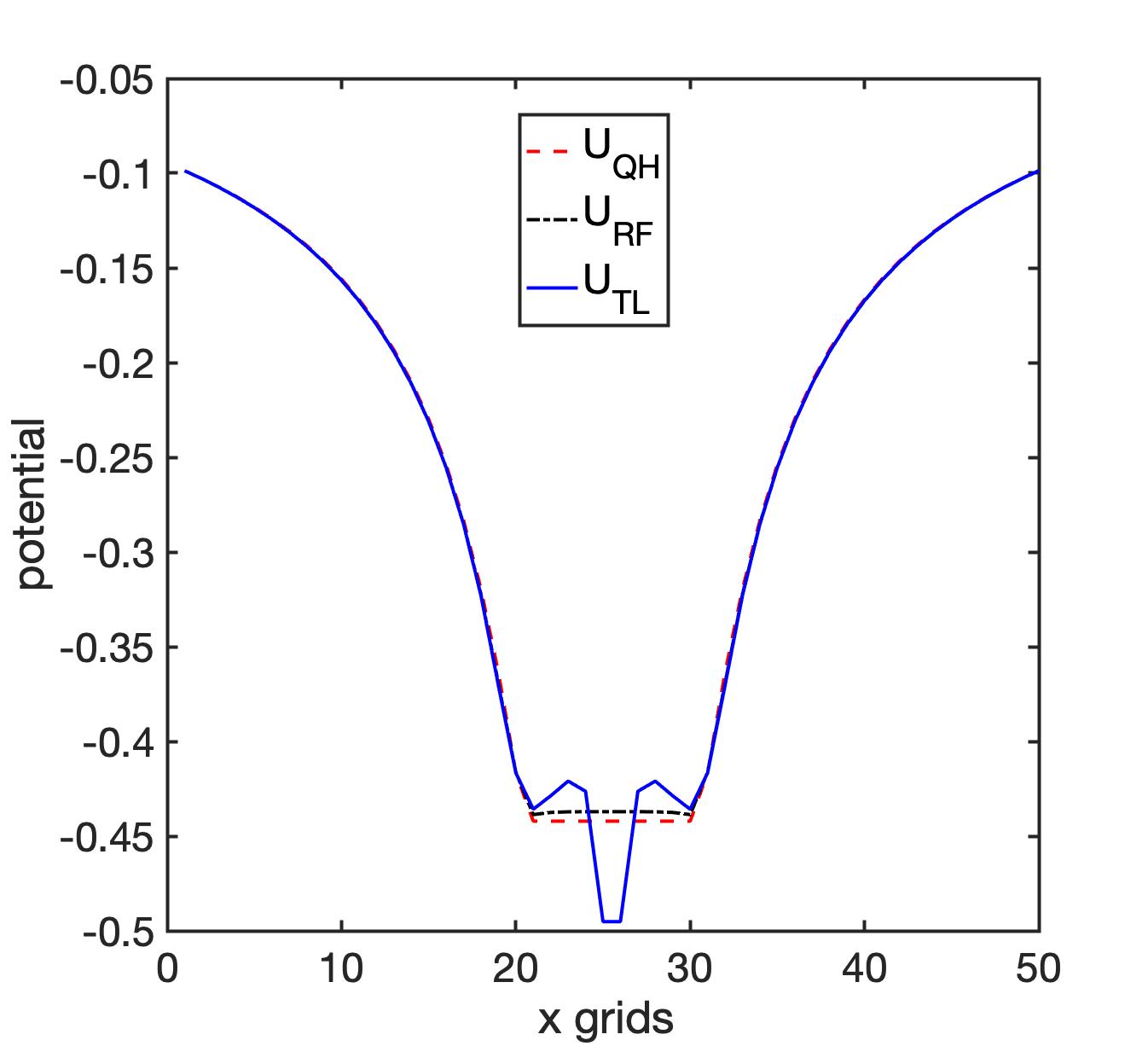}
  \caption{$N=50$}
\end{subfigure}
\begin{subfigure}{.45\textwidth}
  \centering
  \includegraphics[width=0.95\linewidth]{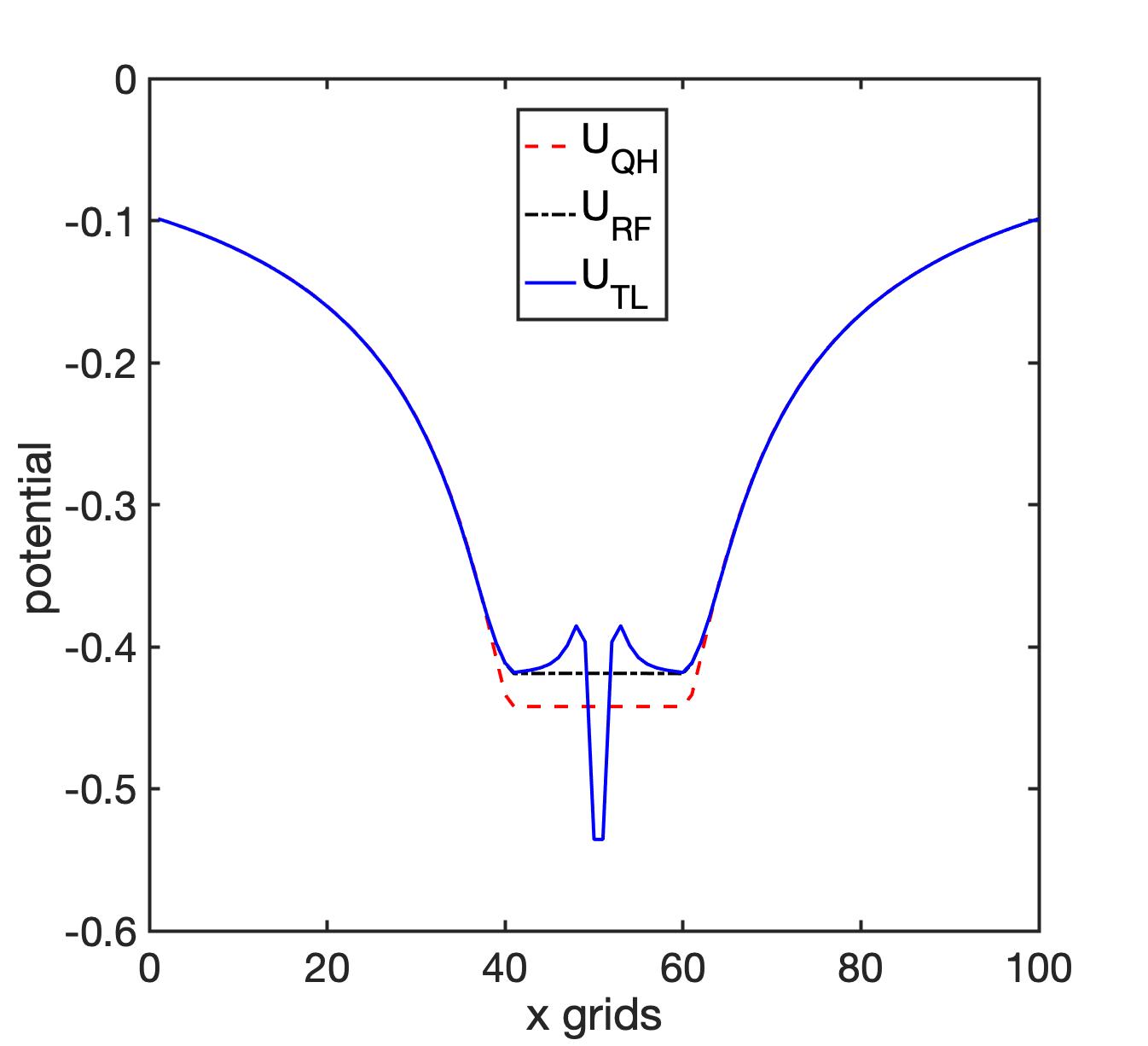}
  \caption{$N=100$}
\end{subfigure} \\
\begin{subfigure}{.45\textwidth}
  \centering
  \includegraphics[width=0.95\linewidth]{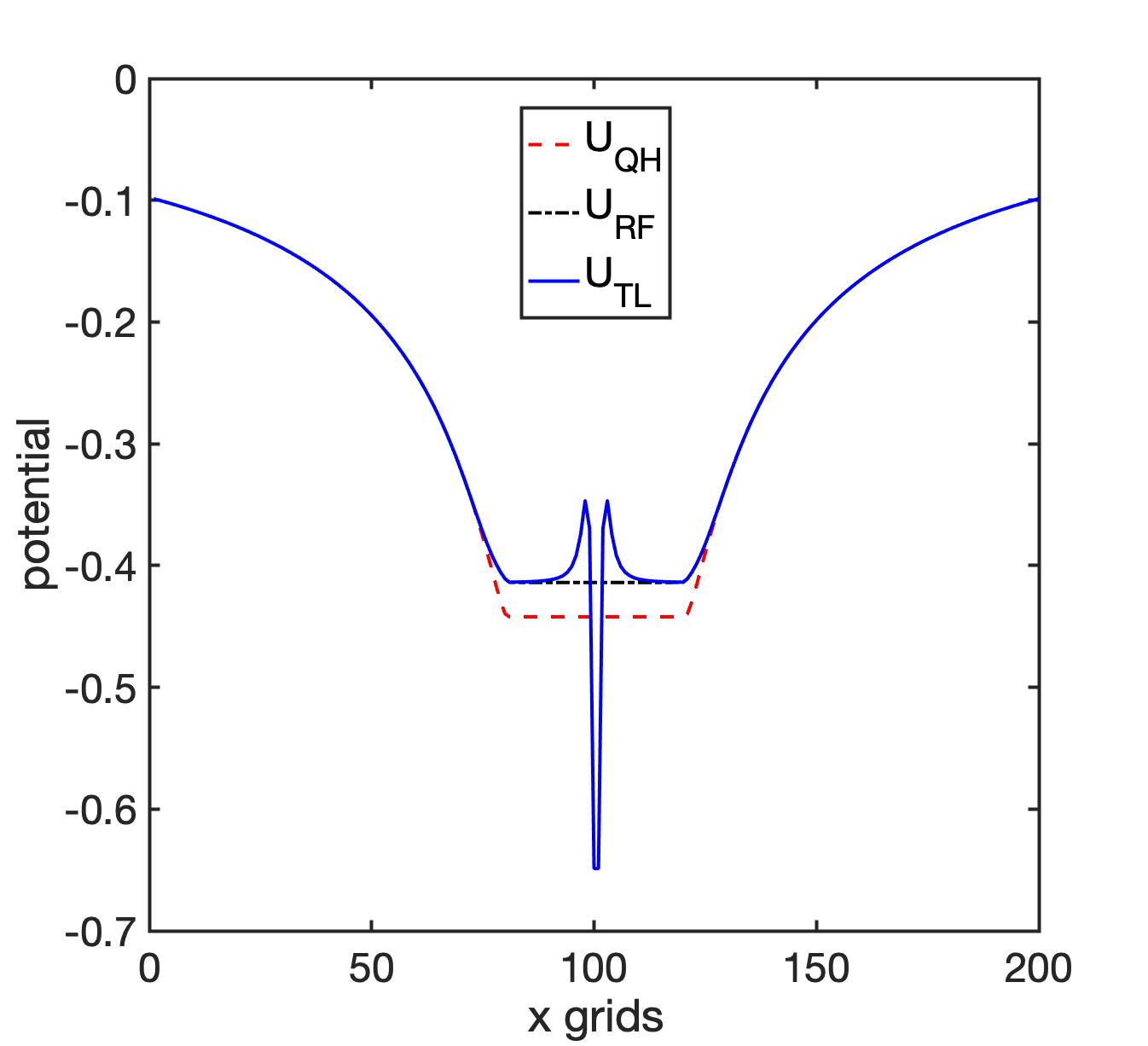}
  \caption{$N=200$}
\end{subfigure}
\begin{subfigure}{.45\textwidth}
  \centering
  \includegraphics[width=0.95\linewidth]{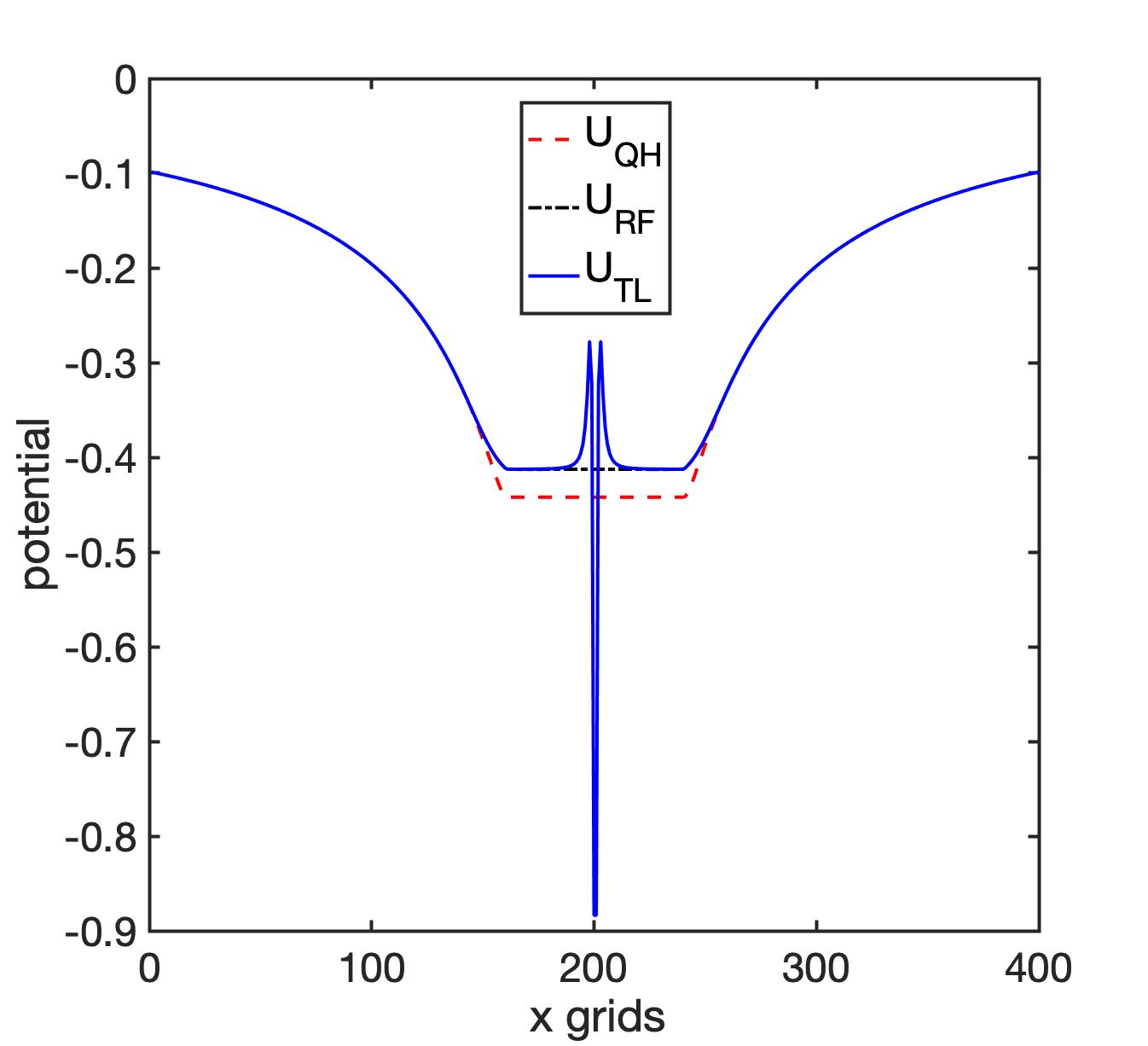}
  \caption{$N=400$}
\end{subfigure}
\caption{Line plots of $u_{QH}$, $u_{RF}$ and $u_{TL}$ along a $x$ line with $y = 0$ and $z = 0$. }
\label{Fig.2d}
\end{figure}

The subtle difference between $u_{RF}$ and $u_{QH}$ lies in the height level of flat potential in $\Omega_i$. To see this, we focus on the difference $|u_{QH} - u_{RF}|$  in Fig. \ref{Fig.3d2} by considering $N=50$, 100, 200, and 400. In all plots, the difference is almost zero outside the sphere $\Omega_i$, and major disagreement only occurs inside the sphere. For a small $N$, like $N=50$, the surface plot of $u_{QH}$ is not completely flat, due to numerical errors. Nevertheless, as $N$ becomes larger, $u_{QH}$ becomes flat enough so that the difference eventually looks like a cylinder with $N=400$. It is noted that the height difference between two solutions approaches to a constant, which is around $0.03$, as shown in Fig. \ref{Fig.3d2} (d). 

The numerical convergence and height difference can also be visualized by depicting three potentials along a $x$ line with $y=0$ and $z=0$, see Fig. \ref{Fig.2d}. For trilinear solution $u_{TL}$, it obviously converges to $u_{RF}$ in most parts, except for near the charge center. However, near the origin, the disagreement between  $u_{TL}$ and $u_{RF}$ increases, suggesting a divergent behavior of trilinear charge distribution. For $u_{QH}$ and $u_{RF}$, their difference becomes negligible away from the sphere $\Omega_i$. Inside the sphere, the height difference indeed approaches to a constant. 

We finally quantitatively compare the difference of three solutions in $L_2$ and $L_{\infty}$ norms for different $N$ in Table \ref{table1}. For the difference between $u_{TL}$ and $u_{RF}$, the $L_2$ norm becomes smaller and smaller. This agrees the above observation that both numerical solutions converge to the same place as $h$ goes to zero. However, the $L_{\infty}$ norm diverges in a rate inverse proportional to $h$, i.e., $O(h^{-1})$. This result fully illustrates how bad the trilinear approximation is. Fortunately, such a difficulty is analytically bypassed in our regularization method. For $u_{QH}$ and $u_{RF}$, we note that the height difference between two solutions inside $\Omega_i$ is actually captured by the $L_{\infty}$ norm, which is 3.04E-2 at $N=400$. In fact, the $L_{\infty}$ norm converges \emph{quadratically} to a constant height difference.  To see this, we take 3.04E-2 as the reference value for the ``exact'' height difference. Then the change in the $L_\infty$ norm is 1.74E-2, 0.42E-2, and 0.09E-2, respectively, for $N=50$, 100, and 200. This obviously is a sequence with $O(h^2)$ convergence, and demonstrates the second order accuracy of the central difference discretization underlying the regularization approach. 

\begin{table}[!tb]
\centering
\caption{The comparison of three solutions' differences.}
\begin{tabular}{|c|c|cc|cc|cc|}
\hline 
&& \multicolumn{2}{c|}{$|u_{QH}-u_{RF}|$} & \multicolumn{2}{c|}{$|u_{QH}-u_{TL}|$} 
 & \multicolumn{2}{c|}{$|u_{RF}-u_{TL}|$} \\
\hline 
$N$ & $h$ & $L_2$ & $L_{\infty}$ & $L_2$ & $L_{\infty}$ & $L_2$ & $L_{\infty}$ \\
\hline
50 & 0.408 & 7.66E-4 & 1.30E-2 & 9.60E-4 & 5.30E-2 & 6.12E-4 & 5.80E-2 \\
\hline
100 & 0.202 & 1.82E-3 & 2.62E-2 & 1.87E-3 & 9.35E-2 & 4.53E-4 & 1.16E-1 \\
\hline
200 & 0.101 &  2.19E-3 & 2.95E-2 & 2.21E-3 & 2.07E-1 & 3.23E-4 & 2.35E-1 \\
\hline
400 & 0.050 &  2.29E-3 & 3.04E-2 & 2.30E-3 & 4.41E-1 & 2.29E-4 & 4.71E-1 \\
\hline
\end{tabular}
\label{table1}
\centering
\end{table}

\section{Conclusion}
A novel regularization approach is introduced for Poisson's equation with singular charge sources and diffuse interfaces, which is the first of its kind in the literature. Through a dual decomposition of potential and dielectric functions, the proposed regularized Poisson equation for the reaction field potential has the same elliptic operator with a smooth source function, which can be easily solved by common numerical methods. For a simple spherical problem, the regularization method is validated by comparing with a semi-analytical method and conventional trilineary distribution method. The further development of the regularization method for the Poisson Boltzmann equation with diffuse interfaces will be reported in the future.

\section*{Acknowledgments}
This research is partially supported by the National Science Foundation (NSF) of USA under grant DMS-1812930.

\end{document}